\documentclass[a4paper,12pt, reqno]{amsart}
\usepackage{a4}

\usepackage{amsfonts,exscale}
\usepackage{mleftright}

\usepackage{srcltx,hyperref}
\usepackage{cite}

\usepackage{color}

\long\def \forget#1{}

\forget{
\pdfoutput=1

\usepackage[pdftex]{hyperref}
\hypersetup{%
pdftitle ={Large sieve inequality},
pdfauthor = {Rajneesh Kumar singh, Stephan Baier},
pdfcreator = {pdfLatex},
pdfproducer={pdfLatex with hyperref}
}
}

\usepackage{amsmath}
\usepackage{amssymb}
\usepackage{amscd}
\usepackage{textcomp}
\usepackage{amsthm}
\theoremstyle{plain}
\newtheorem{Lemma}{Lemma}[section]
\newtheorem{Theorem}[Lemma]{Theorem}

\newtheorem{Corollary}[Lemma]{Corollary}

\theoremstyle{definition}

\newcommand{\DS}{\displaystyle}

\newcommand{\es}{\enspace}

\newcommand{\dpl}{{\mathchoice{\mbox{\rm (\hspace{-0.15em}(}}
{\mbox{\rm (\hspace{-0.15em}(}}
{\mbox{\scriptsize\rm (\hspace{-0.15em}(}}
{\mbox{\tiny\rm (\hspace{-0.15em}(}}}}
\newcommand{\dpr}{{\mathchoice{\mbox{\rm )\hspace{-0.15em})}}
{\mbox{\rm )\hspace{-0.15em})}}
{\mbox{\scriptsize\rm )\hspace{-0.15em})}}
{\mbox{\tiny\rm )\hspace{-0.15em})}}}}
\newcommand{\invlim}[1][]{\ifthenelse{\equal{#1}{}}
{\DS \lim_{\longleftarrow}}
{\DS \lim_{\underset{#1}{\longleftarrow}}}
}
\newcommand{\dirlim}[1][]{\ifthenelse{\equal{#1}{}}
{\DS \lim_{\longrightarrow}}
{\DS \lim_{\underset{#1}{\longrightarrow}}}
}

\newcommand{\ol}[1]{{\overline{#1}}}

\newcommand{\BOne} {{\mathchoice{\hbox{\rm1\kern-2.7pt l\kern.9pt}}
{\hbox{\rm1\kern-2.7pt l\kern.9pt}}
{\hbox{\scriptsize\rm1\kern-2.3pt l\kern.4pt}}
{\hbox{\scriptsize\rm1\kern-2.4pt l\kern.5pt}}}}

\newcommand{\BC}{{\mathbb{C}}}

\newcommand{\BF}{{\mathbb{F}}}

\newcommand{\BR}{{\mathbb{R}}}

\newcommand{\BT}{{\mathbb{T}}}

\newcommand{\sM}{{\mathscr{M}}}

\usepackage{mathrsfs}

\newcommand{\CO}{{\mathcal{O}}}

\newcommand{\Fa}{{\mathfrak{a}}}

\DeclareMathOperator{\Tr}{Tr}
\DeclareMathOperator{\vol}{vol}

\begin{document}
\author{Stephan Baier, Rajneesh Kumar Singh}
\address{Stephan Baier, Ramakrishna Mission Vivekananda Educational Research Institute, G. T. Road, PO Belur Math, Howrah, West Bengal 711202, 
India; email: email$_{-}$baier@yahoo.de}
\address{Rajneesh Kumar Singh, Ramakrishna Mission Vivekananda Educational Research Institute, G. T. Road, PO Belur Math, Howrah, West Bengal 711202, India; rajneeshkumar.s@gmail.com}
\title{Large sieve inequality with power moduli for function fields}
\maketitle
\begin{abstract}
In this paper, we establish a general version of the large sieve with additive characters for restricted sets of moduli 
in arbitrary dimension for function fields. From this, we derive function field versions for the large sieve in high dimensions and for
power moduli. 
\end{abstract}
\tableofcontents
\bigskip

\newpage
\section{Introduction}
The classical large sieve inequality with additive characters asserts that
$$
\sum\limits_{q\le Q} \sum\limits_{\substack{a=1\\ (a,q)=1}}^q \left| \sum\limits_{M<n\le M+N} a_n e\left(n\cdot \frac{a}{q}\right) \right|^2
\le (N-1+Q^2)\sum\limits_{M<n\le M+N} |a_n|^2,
$$
where $Q,N\in \mathbb{N}$ and $M\in \mathbb{Z}$. 
As the name tells us, the large sieve is useful for sieving problems. Indeed, there is an arithmetic form of the large sieve due to Montgomery
\cite{Mon}. Both versions of the large sieve have numerous applications in analytic number theory. The large sieve with resticted sets of 
moduli $q$, in particular power moduli, was considered in a series of papers by Baier, Zhao and Halupczok (see \cite{Bai}, \cite{Zh1}, \cite{BZ1},
\cite{BZ2}, \cite{Hal}), and these results turned out to be useful tools for 
applications (see \cite{BPS} and \cite{BFKS}). In particular, for square moduli, it was established in \cite{BZ2} that 
\begin{equation*}
\begin{split}
& \sum\limits_{q\le Q} \sum\limits_{\substack{a=1\\ (a,q)=1}}^{q^2} \left| \sum\limits_{M<n\le M+N} a_n e\left(n\cdot \frac{a}{q^2}\right) 
\right|^2 \\
\ll_{\varepsilon} & (QN)^{\varepsilon}\left(N+Q^3+\min\left\{N\sqrt{Q},\sqrt{N}Q^2\right\}\right)\cdot \sum\limits_{M<n\le M+N} |a_n|^2.
\end{split}
\end{equation*}
In \cite{Zh1}, Zhao conjectured that the bound
\begin{equation} \label{conjec}
\sum\limits_{q\le Q} \sum\limits_{\substack{a=1\\ (a,q)=1}}^{q^k} \left| \sum\limits_{M<n\le M+N} a_n e\left(n\cdot \frac{a}{q^k}\right) \right|^2
\ll_{\varepsilon} (QN)^{\varepsilon}\left(N+Q^{k+1}\right)\cdot \sum\limits_{M<n\le M+N} |a_n|^2
\end{equation}
should hold for for $k$-th power moduli. This conjecture is still open for every $k$. 
A version of the classical large sieve in higher dimensions was proved by Gallagher \cite{Gal} and subsequently improved by Zhao \cite{Zh2}. The 
last-named author's version states that 
\begin{equation} \label{highdimensionzhao}
\begin{split}
\sum\limits_{\substack{\beta\in \mathbb{R}^n/\mathbb{Z}^n\\ \mbox{\scriptsize \rm ord}(\beta)\le Q}} 
\left|\sum\limits_{\substack{\alpha=(\alpha_1,...,\alpha_n)\in \mathbb{Z}^n\\ \max_{1\le i\le n}|\alpha_i|\le N}}
c(\alpha) \cdot e(\alpha\cdot \beta) \right|^2 \ll & \left(N^n+Q^2N^{n-1}+Q^{n+1}\right) \times\\
& \sum\limits_{\substack{\alpha=(\alpha_1,...,\alpha_n)\in \mathbb{Z}^n\\ \max_{1\le i\le n}|\alpha_i|\le N}}
|c(\alpha)|^2.
\end{split}
\end{equation}
The function field analogue of the arithmetical form of the large sieve was established by Hsu \cite{Hsu} in arbitrary dimension. In this paper,
we prove a large sieve inequality with additive characters in arbitrary dimension with restricted set of moduli for function fields. 
As consequences of this, complete analogues
of \eqref{conjec} and \eqref{highdimensionzhao} for function fields will be established.

\section{Notation} 
The following notations and conventions are used throughout paper:

$\bullet f = O(g)$ means $|f| \leq cg$ for some unspecified positive constant $c$.\\

$\bullet f\ll g$ means $f = O(g)$.\\

$\bullet f\asymp g$ means $c_1g\leq f \leq c_2g$ for some unspecified positive constants $c_1$ and $c_2$. Unless otherwise stated, all implied constants in $\ll, O$ and $\asymp$ are absolute.

\bigskip
Let $\BF_q$ be a fixed finite field with $q$ elements of characteristic $p$ and let $\Tr : \BF_q \to \BF_p$ be the trace map.

\bigskip
\noindent
Let $\BF_q(t)_{\infty}$ be the completion of $\BF_q(t)$ at $\infty$ (i.e. $\BF_q\dpl 1/t\dpr$), let $\CO_{\infty}$ be the 
maximal compact subring of $\BF_q(t)_{\infty}$, and let $\sM_{\infty}$
be the maximal ideal of $\CO_{\infty}$. 
The absolute value $|\cdot |_{\infty}$ of $\BF_q(t)_{\infty}$ is defined by 

\begin{align*}
\mathrel \bigg | \sum_{i= -\infty}^n a_i t^i\bigg| = q^n, \es \text{if} \es 0 \neq a_n \in \BF_q.
\end{align*}

\noindent
The non-trivial additive character $E: \BF_q \to \BC^\times$ is defined by 

\begin{align*}
E(x) = \exp \bigg\{\frac{2\pi i}{p} \Tr(x)\bigg\},
\end{align*}
and the map $e: \BF_q(t)_\infty \to \BC^\times$ is defined by 
\begin{align*}
e\bigg(\sum_{i= -\infty}^n a_i t^i\bigg) = E(a_{-1}).
\end{align*}
This map $e$ is also a non-trivial additive character of $\BF_q(t)_\infty$.

\bigskip
\noindent
Given $f = (f_1, f_2,\cdots, f_n)\in \BF_q(t)^n_\infty$, we define the additive character $\Psi_f:  \BF_q(t)^n_\infty \to \BC^\times$ as 
\begin{align*}
\Psi_f((g_1,g_2,\cdots,g_n))&= e(f_1g_1+f_2g_2+\cdots,f_ng_n)\\
&=\prod_{i=1}^ne(f_ig_i)
\end{align*}
for any $(g_1,g_2,\cdots,g_n)\in \BF_q(t)^n_\infty$.

\medskip
\noindent
Suppose $f = (f_1, f_2,\cdots, f_n)\in \BF_q(t)^n_\infty$, we define metric $|\cdot |_\infty$ on $\BF_q(t)^n_\infty$ as 
\begin{align*}
|f |_\infty = \sup\{|f_1 |_\infty,\ |f_2 |_\infty, \cdots,|f_n |_\infty\}.
\end{align*}
 Given an integer $N\geq 0$, the $N$-ball $B(f,N)$ is defined by 
\begin{align*}
B(f,N) = \{g\in  \BF_q(t)^n_\infty \ | \ |g-f|_\infty \leq q^N  \}.
\end{align*}
We view $\BF_q[t]^n\subset \BF_q(t)^n_\infty$ as a lattice of rank $n$ over $\BF_q[t]$, and define the $n$-torus to be $\BT^n = \BF_q(t)^n_\infty/\BF_q[t]^n $. The metric on $\BT^n$ is given by 

\begin{align*}
\|f\|_\infty = \inf_{f'\sim f}|f'|_\infty, 
\end{align*}
where $f'\sim f$ means that $f'\in f + \BF_q[t]^n$. Note that 
$\BT^n$ is a compact Hausdorff space and for all $f\in \BT^n,\ \|f\|_\infty\leq 1/q$.

\section{Preliminaries}
In this section, we collect the results that we need in the course of this paper.
The next Theorem \ref{Thm1} says that the Pontryagin duality holds for rational function fields.

\begin{Theorem}[Theorem 2.1 in \cite{Hsu}] \label{Thm1} Let $B:\BF_q(t)^n_\infty \times \BF_q(t)^n_\infty \to \BC^\times$ be defined as 
$$
B(f,g) = e(f_1g_1+f_2g_2+\cdots,f_ng_n) = \Psi_f(g),
$$
for $f = (f_1, f_2,\cdots, f_n)$ and $g=(g_1,g_2,\cdots,g_n)\in \BF_q(t)^n_\infty$. The Pontryagin duality $\widehat{\BF_q(t)^n_\infty} =\BF_q(t)^n_\infty$ is induced by $B$. Moreover, $\BF_q[t]^n$ is a discrete subgroup of $\BF_q(t)^n_\infty$, and $\widehat{\BT^n}\cong \{\Psi_f\ | \ f \in \BF_q[t]^n\}\cong \BF_q[t]^n$.
\end{Theorem}

\medskip
\noindent
For the locally compact topological space $\BF_q(t)_\infty$, we normalize the Haar measure so that $\mu(\sM_\infty)=1$. For a given locally constant function  $\varphi : \BF_q(t)^n_\infty \to \BC$ with compact support, the Fourier transform $\hat \varphi$ is defined as usual by 
\begin{align*}
\hat \varphi(f) = \int_{\BF_q(t)^n_\infty}\varphi(g)\ol{\Psi_f(g)}dg \ \text{for any} \ f\in\BF_q(t)^n_\infty.
\end{align*}

Using the above Theorem \ref{Thm1}, we can now state the next  theorem \ref{Thm2} which tells us that Poisson summation formula holds for rational function fields and the proof is 
standard and we shall omit it.
\begin{Theorem}[Poisson Summation Formula] \label{Thm2} Let $\Lambda$ be complete lattice in $\BF_q(t)^n_\infty$ and let 
\begin{align*}
\Lambda'= \Big\{g\in\BF_q(t)^n_\infty \mathrel \Big| B(f,g)=1\ \text{for all}\ f\in\BF_q(t)^n_\infty \Big\}
\end{align*}
be the lattice dual to $\Lambda$.  Let $f: \BF_q(t)^n_\infty \to \BC$ be a function such that 
$$
\sum_{a\in \BF_q(t)}|f(x+a)|
$$ 
is uniformly convergent on compact subsets and 
$$
\sum_{a\in \BF_q(t)}|\hat f(a)|
$$ 
is convergent. Then 
\begin{align*}
\sum_{a\in \Lambda}f(a) = \frac{1}{\vol(\Lambda)}\sum_{a'\in \Lambda'}\hat f(a'),
\end{align*}  
where $\vol(\Lambda)$ is the volume of a fundamental mesh of $\Lambda$.

\end{Theorem}

Next  we quote the duality principle.

\begin{Lemma}[Duality Principle, Theorem 288 in \cite{HLP}] Let $T = [t_{mn}]$ be a finite matrix with complex entries. The following two statements are equivalent: 
\begin{enumerate}
\item For any complex numbers $\{a_n\}$, we have 
\begin{align*}
\sum_m \mathrel \Big |\sum_n a_n t_{mn}\Big |^2 \leq \Delta \sum_n |a_n|^2.
\end{align*}
\item For any complex numbers $\{b_n\}$, we have
\begin{align*}
\sum_n \mathrel \Big |\sum_m b_m t_{mn}\Big |^2 \leq \Delta \sum_m |b_m|^2. 
\end{align*}
\end{enumerate}
\end{Lemma}


\section{Large sieve with additive characters}
In analogy to the classical large sieve with restricted sets of moduli, we are interested in having an estimate  of the following kind:
\begin{align}\label{Ineq1}
\sum_{\substack{G \in S}} \sum_{\substack{\sigma \bmod G,\\ \sigma \ \text{proper}\\ \mbox{\scriptsize ord}(\sigma)\le Q}}\mathrel \bigg |\sum_{g \in B(0, N)\cap \BF_q[t]^n} 
a_g \sigma(g+\Fa)\bigg |^2 \leq \Delta \sum_{g \in B(0, N)\cap \BF_q[t]^n} |a_g|^2.
\end{align}
Here, $S$ is a set of subgroups in $\BF_q[t]^n$,  $\sigma$ is an additive character for $\BF_q[t]^n/G$, proper means not a character for a 
subgroup $H$ with $G\subsetneq H$, 
ord$(\sigma)$
is the order of $\sigma$ and $(a_g)_{g\in \BF_q[t]^n}$ is a sequence of complex numbers.

The subgroups of $\BF_q[t]$ are all principle ideals $(f)$, where we can choose  $f$ monic. Hence, the subgroups in $\BF_q[t]^n$ are of the form
$$
(f_1)\times\cdots\times (f_n),
$$
where $f_1,...,f_n$ are monic polynomials in $\BF_q[t]$. 
The proper additive characters $\sigma$ for $\BF_q[t]^n/\Fa$ take then the values 
\begin{align*}
\sigma(g+\Fa)  = e\Big(g\cdot \left(\frac{r_1}{f_1},...,\frac{r_n}{f_n}\right)\Big), \ \text{for some}\ r_i \ \text{with} \ (r_i,f_i)=1.
\end{align*}
We note that
$$
\mbox{ord}(\sigma)=\mbox{deg}(\mbox{lcm}(f_1,...,f_n)).
$$
We denote by $T$ the  left hand side in inequality \eqref{Ineq1} and write
$$
f=(f_1,...,f_n), \quad r=(r_1,...,r_n),\quad \frac{r}{f}=\left(\frac{r_1}{f_1},...,\frac{r_n}{f_n}\right),
$$
$$
(r,f)=1 \mbox{ if } (r_i,f_i)=1 \mbox{ for all } i=1,...,n
$$
and 
$$
\quad F=\mbox{lcm}(f_1,...,f_n),
$$
where $\mbox{lcm}(f_1,...,f_n)$ is the monic polynomial of smallest degree which is divisible by $f_1,...,f_n$. 
Hence we have 
\begin{align*}
T = \sum_{\substack{f \in \tilde S, \\ \deg  F\leq Q}} \sum_{\substack{r \bmod f,\\ (r,f)=1}}\mathrel \bigg 
|\sum_{g \in B(0, N)\cap \BF_q[t]^n} a_g 
e\Big(g\cdot \frac{r}{f}\Big)\bigg |^2.
\end{align*}
Here $\tilde S$ is the set of all $n$-tuples of monic polynomials $f$ such that
\begin{align*}
S = \big\{(f_1)\times\cdots \times (f_n) \ : \ f\in \tilde S\big\}.
\end{align*}
By duality principle, 
\begin{align*}
T \leq \Delta \sum_{g \in B(0, N)\cap \BF_q[t]^n} |a_g |^2
\end{align*}
for all $(a_g)_{g\in \BF_q[t]^n}$ if and only if 
\begin{align*}
T': = \sum_{g \in B(0, N)\cap \BF_q[t]^n}\mathrel \bigg|\sum_{\substack{f \in \tilde S, \\ \deg  F\leq Q}} 
\sum_{\substack{r \bmod f,\\ (r,f)=1}} b_{f,r}e\Big(g\cdot \frac{r}{f}\Big)\bigg|^2 \leq \Delta \sum_{f,r}|b_{f,r}|^2
\end{align*}
for all $(b_{f,r})_{f\in \tilde{S}, \ r\ \text{with}\ (r,f)=1 }$.

\bigskip
\noindent
\section{A general large sieve inequality}
Now consider, more generally, the sum 
\begin{align*}
\sum_{g \in B(0, N)\cap \BF_q[t]^n}\mathrel \bigg|\sum_{i=1}^{R}b_i \ e(g\cdot X_i)\bigg|^2
\end{align*}
for $X_1,\cdots, X_R \in \BF_q(t)_\infty^n$. Let $\Phi : \BF_q(t)_\infty^n\to \BR_{>0}$ be a function satisfying $\Phi(x)\geq 1$ if $|x|_\infty \leq 1$. Then 
\begin{align*}
\sum_{g \in B(0, N)\cap \BF_q[t]^n}\mathrel \bigg|\sum_{i=1}^{R}b_i \ e(g\cdot X_i)\bigg|^2 &\leq  \sum_{g \in  \BF_q[t]^n}
\Phi\Big(\frac{g}{t^N}\Big)\mathrel \bigg|\sum_{i=1}^{R}b_i \ e(g\cdot X_i)\bigg|^2 \\
& = \sum_{i_1,i_2 =1}^{R}b_{i_1}\ol b_{i_2}\sum_{g \in  \BF_q[t]^n}\Phi\Big(\frac{g}{t^N}\Big)e\Big(g\cdot (X_{i_1}-X_{i_2})\Big).
\end{align*}
Now we denote by 
\[
f(g)= \Phi\Big(\frac{g}{t^N}\Big)\cdot e\Big(g\cdot (X_{i_1}-X_{i_2})\Big),
\]
then we have 
\begin{align*}
\hat{f}(g) = \int_{\BF_q(t)_\infty^n}\Phi\Big(\frac{h}{t^N}\Big)e\Big(h\cdot (X_{i_1}-X_{i_2})\Big) e(-h\cdot g)dh \ 
\text{for any}\ g \in \BF_q(t)_\infty,
\end{align*}
where  $\hat{f}$ is  the Fourier transform of $f$. Replacing $h$ by $t^Ns$
\begin{align*}
\hat{f}(g)  &=  \left|t^N\right|_\infty^n \cdot \int_{\BF_q(t)_\infty^n}\Phi(s) \cdot e\Big(t^Ns\cdot (X_{i_1}-X_{i_2}-g)\Big)ds \\
 &=  q^{nN}\cdot \hat\Phi\Big(-t^N(X_{i_1}-X_{i_2}-g)\Big).
\end{align*}
Now we choose $\Phi$ exactly like 
\begin{align*}
\Phi(x): & = 
\begin{cases}
1, &  \text{if} \ |x|_\infty \leq 1, \\
0, & \text{otherwise}
\end{cases} \\
&=
\begin{cases}
1, &  \text{if} \ |x/t|_\infty \leq 1/q, \\
0, & \text{otherwise}.
\end{cases}
\end{align*}
Then $\Phi(x)= \Phi_1(x/t)$ with 
\[
\Phi_1(x)=
\begin{cases}
1, &  \text{if} \ |x|_\infty \leq 1/q, \\
0, & \text{otherwise}.
\end{cases}
\]
We know from ~\cite[Lemma 2.2]{Hsu} that $\hat{\Phi}_1=\Phi_1$. Hence,
\begin{align*}
\hat \Phi(x) = q^n \hat \Phi_1(tx)= \begin{cases}
q, &  \text{if} \ |x|_\infty \leq q^{-2}, \\
0, & \text{otherwise}.
\end{cases}
\end{align*}
From the calculations above and the Poisson summation formula  it follows that 
\begin{align*}
& \sum_{g \in B(0, N)\cap \BF_q[t]^n}\mathrel \bigg|\sum_{i=1}^{R}b_i \ 
e(g\cdot X_i)\bigg|^2\\ &\leq \sum_{i_1,i_2 =1}^{R}b_{i_1}\ol b_{i_2}\sum_{g \in  \BF_q[t]^n}
\Phi\Big(\frac{g}{t^N}\Big)e\Big(g\cdot (X_{i_1}-X_{i_2})\Big)\\
 = &\sum_{i_1,i_2 =1}^{R}b_{i_1}\ol b_{i_2}\sum_{g \in  \BF_q[t]^n}  q^{nN}\cdot \hat\Phi\Big(-t^N(X_{i_1}-X_{i_2}-g)\Big)\\
 = & \sum_{i_1,i_2 =1}^{R}b_{i_1}\ol b_{i_2}\sum_{g \in  \BF_q[t]^n}  q^{n(N+1)}\cdot \hat\Phi_1\Big(-t^{N+1}(X_{i_1}-X_{i_2}-g)\Big)\\ 
 = & \es  q^{n(N+1)}\sum_{\substack{i_1,i_2 =1\\ \| x_{i_1}-x_{i_2}\|\leq q^{-(N+2)}}^{R}}b_{i_1}\ol b_{i_2}\\  
 \leq &\es  \frac{1}{2}\es q^{n(N+1)}\cdot\sum^{R}_{\substack{i_1,i_2 =1\\ \| x_{i_1}-x_{i_2}\|\leq q^{-(N+2)}}}
 \left(|b_{i_1}|^2+ |b_{i_2}|^2\right)\\   
 \leq& \es q^{n(N+1)}\sum_{i_1 =1}^{R}|b_{i_1}|^2\cdot \#\Big\{1\leq i_2\leq R \mathrel\Big| \|X_{i_1}-X_{i_2}\|\leq   
 q^{-(N+2)}\Big\}\\
=  &\es q^{n(N+1)}\cdot\max_{1\leq i_1\leq R}\#\Big\{1\leq i_2\leq R \mathrel\Big| \|X_{i_1}-X_{i_2}\|\leq  q^{-(N+2)}\Big\}\cdot
\sum_{i=1}^{R}|b_i|^2
\end{align*}
where $\|\cdot\|$ is the induced distance on the torus $\BT = \BF_q(t)_\infty^n/\BF_q[t]^n$.

\section{Case of Farey fractions}
Now, we specify the set of $X_i$'s to be of the form 
\begin{align*}
S_{Q} = \Big\{  r/f \in \BF_q(t)^n\ : \ f\in \tilde{S},\  (r,f)=1,\ \deg r_i< \deg f_i \mbox{ for } i=1,...,n, \ \deg  F\le Q\Big\}.
\end{align*}
Hence $S_{Q}$ consists of analogues of Farey fractions of order $Q$ with a restricted set of denominators. 

\bigskip
\noindent
Define 
\begin{align*}
M(Q, N) = \max_{x\in S_{Q}}\#\Big\{\tilde{x} \in S_{Q} \mathrel\Big| \|\tilde{x}-x\|\leq q^{-N}\Big\}.
\end{align*}
Then combining all of the above results, we obtain
\begin{equation} \label{Tineq}
T\le q^{n(N+1)} M(Q,N+2)\cdot \sum\limits_{g \in B(0, N)\cap \BF_q[t]^n} |a_g|^2.
\end{equation}

Now take $x = r/f$ and $\tilde{x} = \tilde{r}/\tilde{f}$, then 
\begin{align*}
\|\tilde{x}-x\|\leq q^{-N} \Longleftrightarrow |\tilde{r}_if_i-r_i\tilde{f}_i |_\infty\leq q^{\deg f_i +\deg \tilde{f}_i-N} 
\quad \mbox{for all } i=1,...,n. 
\end{align*}
Set 
$$
c_i = \tilde{r}_if_i-r_i\tilde{f}_i\quad \mbox{ for } i=1,...,n
$$
and 
$$
\tilde{F}=\tilde{F}(\tilde{f}):=\mbox{lcm}(\tilde{f}_1,...,\tilde{f}_n).
$$
Then
\begin{align*}
M(Q, N) =& \max_{{x\in S_Q}}\#\Big\{\tilde{x} \in S_Q \mathrel\Big| |c_i|_\infty\leq q^{\deg f_i+\deg \tilde{f}_i -N} 
\mbox{ for } i=1,...,n\Big\}\\
= & \max_{{x\in S_Q}} \sum_{\substack{\tilde{f}\in \tilde S\\ \deg \tilde{F}\leq Q}}\es \prod\limits_{i=1}^n
\left(\sum_{\substack{c_i\in \BF_q[t]\\ |c_i|_\infty\leq q^{\deg f_i+\deg \tilde{f}_i -N}\\ 
c_i\equiv -r_i\tilde{f_i} \bmod f_i}}1\right).
\end{align*}
Now we have 
$$
\DS \sum_{\substack{c_i\in \BF_q[t]\\|c_i|_\infty\leq q^{\deg f_i+\deg \tilde{f}_i -N}\\ c_i\equiv -r_i\tilde{f_i} \bmod f_i}}1
\es  = 
q^{\deg \tilde{f}_i-N}
$$ 
if $\deg \tilde{f}_i\geq N$. But if $\deg \tilde{f}_i< N$, then necessarily $c_i=0$ and hence $r_i= \tilde{r}_i$ and $f_i=\tilde{f}_i$.
It follows that 
\begin{align*}
M(Q, N)\leq &  \max_{f\in \tilde{S}}   
\sum_{\substack{\tilde{f}\in \tilde S\\ \deg \tilde{F}\leq Q}}\es 
\prod\limits_{i=1}^n\Bigg(\delta(f_i,\tilde{f}_i)+q^{\deg \tilde{f_i}-N} \Bigg),
\end{align*}
where 
$$
\delta(f_i,\tilde{f}_i)=\begin{cases} 1 & \mbox{ if } f_i=\tilde{f}_i\\ 0 & \mbox{ if } f_i\not=\tilde{f}_i. \end{cases}
$$
Combining this with \eqref{Tineq}, we obtain the following. 

\begin{Theorem} \label{generalcase} We have 
\begin{align*}
T \leq& \es  q^{n(N+1)} \cdot \max_{f\in \tilde{S}}   
\sum_{\substack{\tilde{f}\in \tilde S\\ \deg \tilde{F}\leq Q}}\es 
\prod\limits_{i=1}^n\Bigg(\delta(f_i,\tilde{f}_i)+q^{\deg \tilde{f_i}-(N+2)} \Bigg)\cdot \sum_{g \in B(0, N)\cap \BF_q[t]^n} |a_g|^2.
\end{align*}
\end{Theorem}

This implies the following for the one-dimensional case.

\begin{Corollary} \label{n=1} If $n=1$, then 
\begin{align*}
T \leq \es  \left(q^{N+1} + \left(\sharp S\right) \cdot q^{Q-1}\right) \cdot  
\sum_{g \in B(0, N)\cap \BF_q[t]^n} |a_g|^2.
\end{align*}
\end{Corollary}

\section{Case of $k$-th power moduli}
Now we focus on the case of $k$-th power moduli, i.e. the case when 
\begin{equation} \label{given}
S:=\left\{\left(f_1^k\right)\times \cdots \times \left(f_n^k\right)\ :\ f=(f_1,...,f_n)\in \BF_q[t]^n \mbox{ monic}\right\},
\end{equation}
where $f=(f_1,...,f_n)\in \BF_q[t]^n$ monic means that all of the polynomials $f_1,...,f_n$ are monic. 
For $k,m,n\in \mathbb{N}$ with $m\ge n$, $f=(f_1,...,f_m)\in \BF_q[t]^m$ monic, $X\ge 0$ and $N>0$, we define
$$
\tilde{M}_{f,n,k}\left(X,N\right):= \sum_{\substack{\tilde{f}\in \BF_q[t]^n \ \mbox{\scriptsize \rm monic}\\ \deg \tilde{F}(\tilde{f})\leq X}}\es 
\prod\limits_{i=1}^n\Bigg(\delta(f_i,\tilde{f}_i)+q^{k\deg \tilde{f_i}-N} \Bigg).
$$ 
Then Theorem \ref{generalcase} implies the following. 

\begin{Corollary} \label{kth} If $S$ is given as in \eqref{given}, then  
\begin{align*}
T \leq& \es  q^{n(N+1)} \cdot \max_{f\in \BF_q[t]^n \ \mbox{\scriptsize \rm monic}} \tilde{M}_{f,n,k}\left(\frac{Q}{k},N+2\right)   
 \sum_{g \in B(0, N)\cap \BF_q[t]^n} |a_g|^2.
\end{align*}
\end{Corollary}

To bound $\tilde{M}_{f,n,k}\left(X,N\right)$, we proceed by recursion over $n$. If $n=1$, we obtain
$$
\tilde{M}_{f,n,k}\left(X,N\right)\le 1+q^{(k+1)X-N}.
$$
If $n>1$, then
\begin{equation*}
\begin{split}
\tilde{M}_{f,n,k}(X,N) = & \sum_{\substack{\tilde{f}\in \BF_q[t]^{n-1} \ \mbox{\scriptsize \rm monic}\\ \deg \tilde{F}(\tilde{f})\leq X}}\es 
\prod\limits_{i=1}^{n-1}\Bigg(\delta(f_i,\tilde{f}_i)+q^{k\deg \tilde{f_i}-N} \Bigg) \times\\ &
\sum\limits_{\substack{\tilde{f}_n\in \BF_q[t] \ \mbox{\scriptsize \rm monic}\\ \deg \mbox{\scriptsize \rm lcm}(\tilde{F}(\tilde{f}),\tilde{f}_n)\leq X}}
\left(\delta(f_n,\tilde{f}_n)+q^{k\deg \tilde{f}_n-N}\right)\\
= & \tilde{M}_{f,n-1,k}\left(X,N\right)+ \sum_{\substack{\tilde{f}\in \BF_q[t]^{n-1} \ \mbox{\scriptsize \rm monic}\\ \deg \tilde{F}(\tilde{f})\leq X}}\es
\prod\limits_{i=1}^{n-1}\Bigg(\delta(f_i,\tilde{f}_i)+q^{k\deg \tilde{f_i}-N} \Bigg) \times\\ &
\sum\limits_{\substack{\tilde{f}_n\in \BF_q[t] \ \mbox{\scriptsize \rm monic}\\ \deg \mbox{\scriptsize \rm lcm}(\tilde{F}(\tilde{f}),\tilde{f}_n)\leq X}}
q^{k\deg \tilde{f}_n-N}.
\end{split}
\end{equation*}
We write 
\begin{equation*}
\begin{split}
\sum\limits_{\substack{\tilde{f}_n\in \BF_q[t] \ \mbox{\scriptsize \rm monic}\\ \deg \mbox{\scriptsize \rm lcm}(\tilde{F}(\tilde{f}),\tilde{f}_n)\leq X}}
q^{k\deg \tilde{f}_n-N} \le & \sum\limits_{g|\tilde{F}(\tilde{f}) \ \mbox{\scriptsize \rm monic}}  
\sum\limits_{\substack{\tilde{f}_n\in \BF_q[t] \ \mbox{\scriptsize \rm monic}\\ g|f_n\\ 
\deg \tilde{F}(\tilde{f})+\deg \tilde{f}_n-\deg g\leq X}}
q^{k\deg \tilde{f}_n-N}\\
= & \sum\limits_{g|\tilde{F}(\tilde{f}) \ \mbox{\scriptsize \rm monic}} \sum\limits_{\substack{f_n^{\ast}\in \BF_q[t] \ \mbox{\scriptsize \rm monic}\\  
\deg \tilde{F}(\tilde{f})+\deg f_n^{\ast}\leq X}}
q^{k(\deg f_n^{\ast}+\deg g)-N}\\
\le & \sum\limits_{g|\tilde{F}(\tilde{f}) \ \mbox{\scriptsize \rm monic}} 
\sum\limits_{\substack{f_n^{\ast}\in \BF_q[t] \ \mbox{\scriptsize \rm monic}\\  
\deg f_n^{\ast}\leq X-\deg \tilde{F}(\tilde{f})}}
q^{k(X-\deg \tilde{F}(\tilde{f})+\deg g)-N}\\
\le & q^{(k+1)(X-\deg \tilde{F}(\tilde{f}))-N} \cdot \sum\limits_{g|\tilde{F}(\tilde{f}) \ \mbox{\scriptsize \rm monic}} q^{\deg g}\\
\le & q^{(k+1)(X-\deg \tilde{F}(\tilde{f}))-N} \cdot (q+1)^{\tilde{F}(\tilde{f})}\\
\le & (q+1)^{(k+1)X-k\deg \tilde{F}(\tilde{f})-N}.
\end{split}
\end{equation*}
It follows that
\begin{equation*}
\begin{split}
\tilde{M}_{f,n,k}(X,N)
\le & \tilde{M}_{f,n-1,k}\left(X,N\right)+ \sum_{\substack{\tilde{f}\in \BF_q[t]^{n-1} \ \mbox{\scriptsize \rm monic}\\ \deg \tilde{F}(\tilde{f})\leq X}}\es
\prod\limits_{i=1}^{n-1}\Bigg(\delta(f_i,\tilde{f}_i)+q^{k\deg \tilde{f_i}-N} \Bigg) \times\\
& (q+1)^{(k+1)X-k\deg \tilde{F}(\tilde{f})-N}\\
\le & \tilde{M}_{f,n-1,k}\left(X,N\right)+ \sum\limits_{0\le j\le X} \sum_{\substack{\tilde{f}\in \BF_q[t]^{n-1} \ \mbox{\scriptsize \rm monic}\\ 
\deg \tilde{F}(\tilde{f})\leq j}}\es
\prod\limits_{i=1}^{n-1}\Bigg(\delta(f_i,\tilde{f}_i)+q^{k\deg \tilde{f_i}-N} \Bigg) \times\\
& (q+1)^{(k+1)X-kj-N}
\end{split}
\end{equation*}
which is the same as 
\begin{equation*}
\tilde{M}_{f,n,k}(X,N)\le \tilde{M}_{f,n-1,k}\left(X,N\right)+ \sum\limits_{0\le j\le X} \tilde{M}_{f,n-1,k}(j,N)\cdot (q+1)^{(k+1)X-kj-N}.
\end{equation*}
Now we have a recursive inequality for $\tilde{M}_{f,n,k}(X,N)$. It is easily checked that this gives rise to the following explicit bound for
$\tilde{M}_{f,n,k}(X,N)$.

\begin{Lemma} We have 
$$
\tilde{M}_{f,n,k}(X,N) \ll_{q,n} 1+(q+1)^{kX+(X-N)}+(q+1)^{kX+n(X-N)}.
$$
\end{Lemma}

So our final estimate for $T$ in the case of power moduli is the following. If $S$ is given as in \eqref{given}, then  
\begin{align*}
T \ll_{q,n} & \es  \left((q+1)^{nN}+ (q+1)^{\frac{k+1}{k}\cdot Q+(n-1)N}+(q+1)^{\frac{k+n}{k}\cdot Q}\right)\cdot 
 \sum_{g \in B(0, N)\cap \BF_q[t]^n} |a_g|^2.
\end{align*}
Changing $Q/k$ into $Q$, we arrive at the following. 

\begin{Theorem} Let $k,n\in \mathbb{N}$. Then
\begin{align*}
& \sum_{\substack{f \in \BF[t]^n\ \mbox{\scriptsize \rm monic}, \\ \deg  F\leq Q}} \sum_{\substack{r \bmod f^k,\\ (r,f)=1}}\mathrel \bigg 
|\sum_{g \in B(0, N)\cap \BF_q[t]^n} a_g 
e\Big(g\cdot \frac{r}{f^k}\Big)\bigg |^2\\
\ll_{q,n} & \left((q+1)^{nN}+ (q+1)^{(k+1)Q+(n-1)N}+(q+1)^{(k+n)Q}\right)\cdot 
 \sum_{g \in B(0, N)\cap \BF_q[t]^n} |a_g|^2,
\end{align*}
where 
$$
f^k:=\left(f_1^k,...,f_n^k\right).
$$
\end{Theorem}

In particular, in the case $k=1$ of full moduli, we get the following result on the large sieve for function fields in dimension $n$.

\begin{Corollary} \label{full} Let $n\in \mathbb{N}$. Then
\begin{align*}
& \sum_{\substack{f \in \BF[t]^n\ \mbox{\scriptsize \rm monic}, \\ \deg  F\leq Q}} \sum_{\substack{r \bmod f,\\ (r,f)=1}}\mathrel \bigg 
|\sum_{g \in B(0, N)\cap \BF_q[t]^n} a_g 
e\Big(g\cdot \frac{r}{f}\Big)\bigg |^2\\
\ll_{q,n} & \left((q+1)^{nN}+ (q+1)^{2Q+(n-1)N}+(q+1)^{(n+1)Q}\right)\cdot 
 \sum_{g \in B(0, N)\cap \BF_q[t]^n} |a_g|^2.
\end{align*}
\end{Corollary}

Moreover, in the case $n=1$ of dimension 1, we get the following result on the large sieve for function fields with power moduli.

\begin{Corollary} \label{dim1} Let $k\in \mathbb{N}$. Then
\begin{align*}
& \sum_{\substack{f \in \BF[t]\ \mbox{\scriptsize \rm monic}, \\ \deg f \leq Q}} \sum_{\substack{r \bmod f^k,\\ (r,f)=1}}\mathrel \bigg 
|\sum_{g \in B(0, N)\cap \BF_q[t]} a_g 
e\Big(g\cdot \frac{r}{f^k}\Big)\bigg |^2\\
\ll_{q,n} & \left((q+1)^{N}+ (q+1)^{(k+1)Q}\right)\cdot 
 \sum_{g \in B(0, N)\cap \BF_q[t]} |a_g|^2.
\end{align*}
\end{Corollary}

\section{Comparison with classical large sieve}
Philosophically, the terms $(q+1)^N$ and $(q+1)^Q$ in the large sieve inqualities for function fields above play the same rules as the terms $N$ and $Q$
in the classical large sieve inequalities, respectively. So the potential analogues of the estimates in Corollaries \ref{full} and 
\ref{dim1} in the classical setting are
\begin{equation} \label{highdimensionconj} 
\begin{split}
\sum\limits_{\substack{\beta\in \mathbb{R}^n/\mathbb{Z}^n\\ \mbox{\scriptsize \rm ord}(\beta)\le Q}} 
\left|\sum\limits_{\substack{\alpha=(\alpha_1,...,\alpha_n)\in \mathbb{Z}^n\\ \max_{1\le i\le n}|\alpha_i|\le N}}
c(\alpha) \cdot e(\alpha\cdot \beta) \right|^2 \ll & \left(N^n+Q^2N^{n-1}+Q^{n+1}\right) \times\\ &
\sum\limits_{\substack{\alpha=(\alpha_1,...,\alpha_n)\in \mathbb{Z}^n\\ \max_{1\le i\le n}|\alpha_i|\le N}}
|c(\alpha)|^2
\end{split}
\end{equation}
and
\begin{equation} \label{kthpowerconj}
\sum\limits_{q\le Q} \sum\limits_{\substack{a=1\\ (a,q)=1}}^{q^k} \left|\sum\limits_{M<n\le M+N} a_ne\left(n\cdot\frac{a}{q^k}\right)\right|^2 
\ll \left(N+Q^{k+1}\right)\sum\limits_{M<n\le M+N} |a_n|^2. 
\end{equation}
Indeed, inequality \eqref{highdimensionconj} is precisely \eqref{highdimensionzhao}, established by Zhao,
and \eqref{kthpowerconj} matches \eqref{conjec} with the term $(QN)^{\varepsilon}$ omitted.

\end{document}